%% This is an Amstex file.  Accepted by
%% Acta Arithmetica August 2008, then withdrawn.
%% This is the new version of July 10, 2014, July 11.

\input amstex.tex
\documentstyle{amsppt}
\magnification1200
\hsize=12.5cm
\vsize=18cm
\hoffset=1cm
\voffset=2cm

\def\DJ{\leavevmode\setbox0=\hbox{D}\kern0pt\rlap
{\kern.04em\raise.188\ht0\hbox{-}}D}
\def\dj{\leavevmode
 \setbox0=\hbox{d}\kern0pt\rlap{\kern.215em\raise.46\ht0\hbox{-}}d}

\def\txt#1{{\textstyle{#1}}}
\baselineskip=13pt
\def\hf{{\textstyle{1\over2}}}

\def\d{{\,\roman d}}
\def\e{\varepsilon}
\def\f{\varphi}
\def\G{\Gamma}

\def\s{\sigma}
\def\t{\theta}
\def\le{\leqslant}
\def\ge{\geqslant}
\def\={\;=\;}
\def\zx{\zeta(\hf+ix)}
\def\zt{\zeta(\hf+it)}
\def\zy{\zeta(\hf+iy)}

\def\D{\Delta}
\def\no{\noindent}
\def\R{\Re{\roman e}\,} \def\I{\Im{\roman m}\,}
\def\z{\zeta}

 \def\t{\theta}
\def\hf{{\textstyle{1\over2}}}
\def\txt#1{{\textstyle{#1}}}
\def\f{\varphi}
\def\Z{{\Cal Z}}
%%%%%%%%%%% Fonts macros %%%%%%%%%%%%
\font\tenmsb=msbm10
\font\sevenmsb=msbm7
\font\fivemsb=msbm5
\newfam\msbfam
\textfont\msbfam=\tenmsb
\scriptfont\msbfam=\sevenmsb
\scriptscriptfont\msbfam=\fivemsb
\def\Bbb#1{{\fam\msbfam #1}}

\def \NN {\Bbb N}

\def \RR {\Bbb R}

\font\ff=cmr8
\def\txt#1{{\textstyle{#1}}}
\baselineskip=13pt

\font\teneufm=eufm10
\font\seveneufm=eufm7
\font\fiveeufm=eufm5
\newfam\eufmfam
\textfont\eufmfam=\teneufm
\scriptfont\eufmfam=\seveneufm
\scriptscriptfont\eufmfam=\fiveeufm
\def\mathfrak#1{{\fam\eufmfam\relax#1}}

\font\tenmsb=msbm10
\font\sevenmsb=msbm7
\font\fivemsb=msbm5
\newfam\msbfam
     \textfont\msbfam=\tenmsb
      \scriptfont\msbfam=\sevenmsb
      \scriptscriptfont\msbfam=\fivemsb
\def\Bbb#1{{\fam\msbfam #1}}

\def \NN {\Bbb N}

\def \RR {\Bbb R}

  \def\rightheadline{{\hfil{\ff
Hybrid moments of the Riemann zeta-function}\hfil\tenrm\folio}}

  \def\leftheadline{{\tenrm\folio\hfil{\ff
   Aleksandar Ivi\'c }\hfil}}
  \def\emptyheadline{\hfil}
  \headline{\ifnum\pageno=1 \emptyheadline\else
  \ifodd\pageno \rightheadline \else \leftheadline\fi\fi}

\font\ff=cmr8
\font\teneufm=eufm10
\font\seveneufm=eufm7
\font\fiveeufm=eufm5
\newfam\eufmfam
\textfont\eufmfam=\teneufm
\scriptfont\eufmfam=\seveneufm
\scriptscriptfont\eufmfam=\fiveeufm
\def\mathfrak#1{{\fam\eufmfam\relax#1}}

\font\tenmsb=msbm10
\font\sevenmsb=msbm7
\font\fivemsb=msbm5
\newfam\msbfam
\textfont\msbfam=\tenmsb
\scriptfont\msbfam=\sevenmsb
\scriptscriptfont\msbfam=\fivemsb
\def\Bbb#1{{\fam\msbfam #1}}

\def \NN {\Bbb N}

\def \RR {\Bbb R}

\def\D{\Delta}

 \def\e{\varepsilon}
\def\no{\noindent} \def\d{\,{\roman d}}
%\topglue1cm
\topmatter
\title
Hybrid moments of the Riemann zeta-function
\endtitle
\author
Aleksandar Ivi\'c
\endauthor
\address
Katedra Matematike RGF-a, Universitet u Beogradu,  \DJ u\v sina 7,
11000 Beograd, Serbia.
\endaddress
\keywords The Riemann zeta-function, power moments, asymptotic formulas,
upper bounds
\endkeywords
\subjclass 11 M 06
\endsubjclass
\email
{\tt  aleksandar.ivic\@rgf.bg.ac.rs, aivic\@matf.bg.ac.rs}
\endemail
\dedicatory
\enddedicatory
\abstract
The ``hybrid" moments
$$
\int_T^{2T}|\zt|^k{\left(\int_{t-G}^{t+G}|\zx|^\ell\d x\right)}^m\d t
\quad\Bigl(T^\e \ll G = G(T) \ll T\Bigr)
$$
of the Riemann zeta-function $\z(s)$ on the critical line $\R s = \hf$
are studied. The expected upper bound for the above expression is
$O_\e(T^{1+\e}G^m)$. This is shown to be true for certain
specific values of $k,\ell,m\in \NN$, and the explicitly determined
range of $G = G(T;k,\ell,m)$. The application to a mean square bound for
the Mellin transform function of $|\zx|^4$ is given.
\endabstract
\endtopmatter
\document
\head
1. Introduction
\endhead
Power moments represent one of the most important parts of the
theory of the Riemann zeta-function $\z(s)$, defined as
$$
\z(s) = \sum_{n=1}^\infty n^{-s}\qquad\qquad(\s = \R s > 1),
$$
and otherwise by analytic continuation.
 Of particular significance are the moments on the
``critical line" $\s = \hf$, and a large literature exists on this
subject (see e.g., the monographs [6], [7], [22], [23] and [25]).
Let us define
$$
I_k(T) = \int_0^T|\zt|^{2k}\d t,\leqno(1.1)
$$
where $k\in\RR$ is a fixed, positive number. Naturally one would want to
find an asymptotic formula for $I_k(T)$ for a given $k$, but this is an
extremely difficult problem. Except when $k = 1$ and $k = 2$, no asymptotic
formula for $I_k(T)$ is known yet, although there are plausible conjectures
for such formulas (see e.g., the work of B. Conrey et al. [2]). In the
absence of asymptotic formulas for $I_k(T)$, one would like then to
obtain good upper bounds for $I_k(T)$. A simple bound
for $|\zt|^k$ is (see [7, Theorem 1.2] and [23])
$$
|\zt|^k \ll \log t\int_{t-1}^{t+1}|\zx|^k \d x + 1,\leqno(1.2)
$$
where $k\in\NN$ is fixed. The use of (1.2) allows one to replace a power
of $|\zt|$ by its integral over a suitable (short)
interval. In employing this procedure one obviously loses something,
but on the other hand one gains flexibility from the fact that
explicit upper bound for $I_k(T+G) - I_k(T-G)$ are known only in the case when
$k = 1$ (see Lemma 1) and $k=2$ (see [7, Theorem 5.2] and [22]).
In this way bounds for
$I_{k+m\ell}(T)$ are reduced to the so-called ``hybrid" moments of the type
$$
\int_T^{2T}|\zt|^k{\left(\int_{t-G}^{t+G}|\zx|^\ell\d x\right)}^m\d t
\qquad(k,\ell,m\in\NN),\leqno(1.3)
$$
where $k,\ell,m$ are assumed to be fixed, and $1\ll G = G(T) \ll T$.
The expected bound for the
expression in (1.3) (this is consistent with the hitherto
unproved Lindel\"of hypothesis that $\zt \ll_\e |t|^\e$) is clearly
$$
O_\e(T^{1+\e}G^m).\leqno(1.4)
$$
Here and later $\e \;(>0)$ denotes arbitrarily small constants,
not necessarily the same ones at each occurrence, and $a = O_\e(b)$
(same as $a \ll_\e b$) means that the implied constant depends only on $\e$.
The problem is to find, for given $k,\ell,m$, the range of
$$
G = G(T;k,\ell,m)
$$
for which the integral (1.3) is bounded by (1.4),
and naturally one would like the lower
bound for $G$ to be as small as possible. Note that from general
results (e.g., see K. Ramachandra's monograph [23]) one obtains
that the expression in (1.3) is, for $\;\log\log T \ll G \ll T$,
$$
\gg\; G^m(\log T)^{\ell^2m/4}\int_T^{2T}|\zt|^k\d t \;\gg\; TG^m(\log
T)^{(\ell^2m+k^2)/4}.\leqno(1.5)
$$
This shows that, up to `$\e$', the bound in (1.4) is indeed best possible.
The (less difficult) case $k=0$ in (1.3) was
investigated by the author in [11] ($\ell = 4$) and [13] ($\ell=2$).
In particular, the former work contains a proof of the bound

\medskip
$$
\int_T^{2T}J_2^m(t,G)\d t \;\ll_\e\; T^{1+\e}\leqno(1.6)
$$
for $T^{1/2+\e} \le G \le T$ if $m = 1,2$; for
$T^{4/7+\e} \le G \le T$ if
$m=3$, and for $T^{3/5+\e} \le G \le T$ if $m=4$, where
$$
J_k(T,G) := {1\over \sqrt{\pi}G}\int_{-\infty}^\infty
|\z(\hf + iT + iu)|^{2k}{\roman e}^{-(u/G)^2}\d u\quad(k>0, 1 \ll G \ll T).
\leqno(1.7)
$$
The bound (1.6) in the above range  was obtained
in [11] by employing Y. Motohashi's
explicit formula (e.g., see [7] and [22]) for $J_2(t,G)$, which contains
quantities from the spectral theory of the non-Euclidean Laplacian.

\smallskip
As for the applications of bounds for (1.3), note that the case
(this is $k = \ell = 4, m=1$) of the hybrid integral
$$
\int_T^{2T}|\zt|^4\int_{t-G}^{t+G}|\zx|^4\d x\d t
\leqno(1.8)
$$
appeared in [19] in connection with mean square bounds for the
Mellin transform function, defined initially by
$$
{\Cal Z}_2(s) = \int_1^\infty |\zx|^4x^{-s}\d x\qquad(\Re s = \s >1),\leqno(1.9)
$$
and otherwise by analytic continuation. The functions ${\Cal Z}_k(s)$
(in the general case $|\zx|^4$ is replaced by $|\zx|^{2k}$ for
 $\s > \s(k) \,(>1)$ with suitable $\s(k)$) are of great
 importance in the theory of power moments
of $\zt$ (see e.g., [9],  [19]). It was shown by the author in [12] that
$$
\int_1^{T}|\Z_2(\s+it)|^2\d t \;\ll_\e\; T^{{15-12\s\over5}+\e}
\quad\qquad(\txt{5\over6} \le \s \le \txt{5\over4}),\leqno(1.10)
$$
which is the sharpest bound for the range in question.

\smallskip

We shall obtain results on the integral in (1.3)
when $k,\ell$ equal 2 or 4, which is logical, since it is in these cases
that we have good information on $I_k(T)$. Namely let, for $k\in\NN$ fixed,
$$
I_k(T) = \int_0^T\vert\zeta(\txt{1\over 2} + it)\vert^{2k}\d t =
T\,P_{k^2}(\log T) + E_k(T), \leqno(1.11)
$$
where for some suitable coefficients $a_{j,k}$ one has
$$
P_{k^2}(y) = \sum_{j=0}^{k^2}a_{j,k}y^j,  \leqno(1.12)
$$
and $E_k(T)$ is to be considered as the error term in (1.11).
An extensive literature exists on $E_k(T)$, especially on
$E_1(T) \equiv E(T)$ (see F.V. Atkinson's classical paper [1]),
and the reader is referred to [7] for a
comprehensive account. It is known that
($\gamma = -\G'(1) = 0.5772157\ldots$
is Euler's constant)
$$
P_1(y) = y + 2\gamma - 1 - \log(2\pi),
$$
and $P_4(y)$ is a quartic polynomial in $y$ whose leading coefficient
equals $1/(2\pi^2)$. This was  obtained in A.E. Ingham's classical work [5].
For  an explicit evaluation of all the
coefficients of $P_4(y)$ see e.g., the author's paper [8]. One hopes that
$$
E_k(T) \= o(T) \qquad (T \to \infty) \leqno(1.13)
$$
will hold for each fixed integer $k \ge 1$, which implies the
Lindel\"of hypothesis that $\zt \ll_\e |t|^\e$. So far (1.13) is
known to be true only in the cases $k = 1$ and $k = 2$, when
$E_k(T)$ is a true error term in the asymptotic formula (1.11). In
particular we have (see [6],[7], [16], [17], [22] and [24]) $E(T) \ll_\e
T^{\t+\e}$ for some $\t$ satisfying ${1\over4} \le \t < {1\over3}$, and
$E(T) = \Omega_\pm(T^{1/4})$. We also have $E_2(T) = \Omega_\pm(\sqrt{T})$
and the bounds (op. cit.)
$$
E_2(T) \ll T^{2/3}\log^8T,\quad \int_0^TE_2^2(t)\d t \ll T^2\log^{22}T.
\leqno(1.14)
$$
As usual, $f(x) = \Omega_\pm(g(x))$ for a given $g(x)\;(>0$ for $x > x_0 > 0$)
means that
$$
\limsup_{x\to\infty}f(x)/g(x)>0,\qquad \liminf_{x\to\infty}f(x)/g(x)<0.
$$

%\bigskip
\head
2. Statement of results
\endhead
Before we state explicitly our results note that we have the bounds
$$
\int_{T-G}^{T+G}|\zt|^2\d t \;\ll\; G\log T\qquad(T^{1/3} \ll G = G(T) \ll T),
\leqno(2.1)
$$
and
$$
\int_{T-G}^{T+G}|\zt|^4\d t \;\ll_\e\; GT^\e\qquad(T^{2/3} \ll G = G(T) \ll T).
\leqno(2.2)
$$
This easily follows from estimates on $E(T)$ and $E_2(T)$ mentioned at the
end of the last section. It means that we can restrict ourselves to
the range $G \ll T^{1/3}$ when $\ell=2$ in (1.3),
and to the range $G \ll T^{2/3}$ when $\ell=4$. This will
be implicitly assumed in the proofs of our results, which are contained
in

\medskip
THEOREM 1. {\it We have}
$$
\int\limits_T^{2T}|\zt|^2\int\limits_{t-G}^{t+G}|\zx|^2\d x\,\d t \ll TG\log^2T
\quad\Bigl(T^\e \ll G = G(T) \ll T\Bigr),\leqno(2.3)
$$
$$
\int\limits_T^{2T}|\zt|^2\int\limits_{t-G}^{t+G}|\zx|^4\d x\,\d t
\ll_\e T^{1+\e}G
\quad(T^{\frac{3}{10}} \ll G = G(T) \ll T),\leqno(2.4)
$$
$$
\int\limits_T^{2T}|\zt|^2{\left(\int_{t-G}^{t+G}|\zx|^2\d x\right)}^2
\,\d t \ll_\e T^{1+\e}G^2
\;(T^{{1\over7}+\e} \ll G = G(T) \ll T),\leqno(2.5)
$$
$$
\int\limits_T^{2T}|\zt|^2\Biggl(\;\int\limits_{t-G}^{t+G}|\zx|^2\d x\Biggr)^3
\,\d t \ll_\e T^{1+\e}G^3
\quad(T^{{1\over5}+\e} \ll G = G(T) \ll T),\leqno(2.6)
$$
{\it and for $\;T^{7\over12}\log^CT \ll G = G(T) \ll T\;$ we have}
$$
\int\limits_T^{2T}|\zt|^2\Biggl(\int_{t-G}^{t+G}|\zx|^4\d x\Biggr)^2
\,\d t \;\ll\; TG^2\log^{9}T.\leqno(2.7)
$$

\medskip
THEOREM 2. {\it We have, for $1 \ll G = G(T) \ll T$ and some $C>0$,}
$$
\int\limits_T^{2T}|\zt|^4\int\limits_{t-G}^{t+G}|\zx|^4\d x\,\d t \ll
\log^CT\Bigl(TG + \min(T^{5/3},T^2G^{-1})\Bigr),
\leqno(2.8)
$$
$$
\int\limits_T^{2T}|\zt|^4{\left(\int_{t-G}^{t+G}|\zx|^2\d x\right)}^2\,\d t
\ll_\e T^{1+\e}G^2\;(T^{7\over24} \le G = G(T) \ll T),\leqno(2.9)
$$
$$
\int\limits_T^{2T}|\zt|^4{\left(\int_{t-G}^{t+G}|\zx|^4\d x\right)}^2\,\d t
\ll_\e T^{1+\e}G^2\;(T^{5\over9} \le G = G(T) \ll T).\leqno(2.10)
$$

\medskip
To assess the strength of our results note, for example, that (1.2)
and (2.8) give
$$\eqalign{
\int_T^{2T}|\zt|^8\d t &\ll \log T\int_T^{2T}|\zt|^4
\int_{t-G}^{t+G}|\zx|^4\d x\,\d t  \cr&+ \int_T^{2T}|\zt|^4\d t
\ll T^{3/2}\log^CT\qquad(G = T^{1/2}).\cr}\leqno(2.11)
$$
The bound in (2.11), which follows easily by the Cauchy-Schwarz
inequality for integrals from estimates of the
fourth and twelfth moment of $|\zt|$
(see [4] and [6, Chapter 8]), is the strongest known bound for the eighth
moment of $|\zt|$.

\medskip
Our last result concerns an improvement of (1.10). Let $\rho$
be such a constant for which
$$
\int_0^T|\zt|^8\d t \;\ll_\e\; T^{\rho+\e}\leqno(2.12)
$$
holds. At present we have $1\le\rho\le3/2$. The lower bound
follows from general principles (see [6, Chapter 9]).
The upper bound is a consequence of (2.11), and its improvements
would be very significant. We shall prove, using (2.4), the following

\medskip
THEOREM 3. {\it If ${\Cal Z}_2(s)$ is defined by}
(1.9) {\it and $\rho$ is defined by} (2.12), {\it then}
$$
\int_1^{T}|\Z_2(\s+it)|^2\d t \;\ll_\e\; T^{{4\rho+4-8\s\over3\rho-1}+\e}
\qquad\Bigl({5+\rho\over8} \le \s \le {1+\rho\over2}\Bigr).\leqno(2.13)
$$

\medskip
{\bf Corollary}. We have
$$\eqalign{&
\int_1^{T}|\Z_2(\txt{13\over16}+it)|^2\d t \;\ll_\e\; T^{1+\e}
\quad \Bigl({13\over16}= 0.8125\,\Bigr),\cr&
\int_1^{T}|\Z_2(1+it)|^2\d t \;\ll_\e\; T^{4/7+\e}
\quad\Bigl({4\over7} = 0.571428\ldots\Bigr).\cr}\leqno(2.14)
$$

\medskip
Note that (1.10) gives
$$
\int_1^{T}|\Z_2(\txt{5\over6}+it)|^2\d t \ll_\e T^{1+\e},\quad
\int_1^{T}|\Z_2(1+it)|^2\d t \ll_\e T^{3/5+\e},
$$
while (2.14) improves both of these bounds,
since $13/16 < 5/6$ and $4/7 < 3/5$.

\head
3. The necessary lemmas
\endhead
In this section we shall state some lemmas that are necessary
for the proofs of our theorems. The first is an explicit formula for an
integral involving $|\zt|^2$.

\bigskip
LEMMA 1. {\it For $T^\e \le G = G(T) \le T^{1-\e}$ we have}
$$
\eqalign{&
{1\over\sqrt{\pi}G}\int_{-\infty}^\infty |\z(\hf+iT+iy)|^2
{\roman e}^{-(y/G)^2}\d y = O(\log T) \;+ \cr&
+ \sqrt{2}\sum_{n=1}^\infty (-1)^nd(n)n^{-1/2}
\left({\Bigl({T\over2\pi n} + {1\over4}\Bigr)}^{1/2}
- {1\over2}\right)^{-1/2}\times\cr&\times\exp\left(-G^2\Bigl(
{\roman {arsinh}}\,\sqrt{\pi n\over2T}\,\Bigr)^{2}\right)\sin f(T,n),\cr}
\leqno(3.1)
$$
{\it where $d(n)$ is the number of divisors of $n$, $\roman{ar\,sinh}\,z =
\log(z+\sqrt{z^2+1}\,)$, and
$$
\eqalign{\cr&
f(T,n) = 2T{\roman {arsinh}}\,\Bigl(\sqrt{{\pi n\over 2T}}\,\Bigr) + \sqrt{2\pi nT
+ \pi^2n^2} - \txt{1\over4}\pi \cr&
=  -\txt{1\over4}\pi + 2\sqrt{2\pi nT} +
\txt{1\over6}\sqrt{2\pi^3}n^{3/2}T^{-1/2} + a_5n^{5/2}T^{-3/2} +
a_7n^{7/2}T^{-5/2} + \ldots\,\cr}\leqno(3.2)
$$
for $\,1 \le n\ll T$, where $\;a_{2m-1}$ are suitable constants.}
\medskip
{\bf Proof of Lemma 1}.  The proof of (3.1) (see also [11]) is based on
Y. Motohashi's exact formula [22, Theorem 4.1]. It states that
$$
\eqalign{&
\int_{-\infty}^\infty |\zt|^2g(t)\d t = \int_{-\infty}^\infty
\left[\R\Bigl\{{\G'\over\G}\bigl(\hf+it\bigr)\Bigr\} + 2\gamma - \log(2\pi)\right]
g(t)\d t \cr&+ 2\pi\R(g(\hf i))
+ 4\sum_{n=1}^\infty d(n)\int_0^\infty (y(y+1))^{-1/2}g_c(\log(1+1/y))
\cos(2\pi ny)\d y,\cr}\leqno(3.3)
$$
where
$$
g_c(x) := \int_{-\infty}^\infty g(t)\cos(xt)\d t
$$
is the cosine Fourier transform of $g(t)$. One requires the function $g(r)$
in (3.3) to be real-valued for $r\in\RR$, and that there exists a
large constant $A>0$ such that $g(r)$ is regular and $\ll (|r|+1)^{-A}$ for
$|\I r| \le A$. The choice
$$
g(t) = {1\over\sqrt{\pi} G}{\roman e}^{-(T-t)^2/G^2},
\quad g_c(x) = {\roman e}^{-{1\over4}(Gx)^2}\cos(Tx)
$$
is permissible, and then the integral on the left-hand side of (3.3)
becomes (see (1.7)) $J_1(T,G)$.
The first integral on the right-hand side of
(3.3) is $O(\log T)$, and the second one is evaluated by the
saddle-point method (see e.g., [6, Chapter 2]). A convenient result to
use is [6, Theorem 2.2 and Lemma 15.1], due originally to Atkinson [1]
for the evaluation of exponential integrals $\int_a^b\f(x)\exp(2\pi iF(x))\d x$.
 In the latter only the exponential
factor $\exp(-{1\over4}G^2\log(1+1/y))$ is missing. In the notation
of [1] and [6] we have that the saddle point $x_0$ (root of $F'(x)=0$) satisfies
$$
x_0 = U - {1\over2} = \left({T\over2\pi n} +
{1\over4}\right)^{1/2} - {1\over2},
$$
and the presence of the above exponential factor makes it possible to
truncate the series in (3.3) at $n = TG^{-2}\log T$ with a negligible error.
Furthermore, in the remaining range for $n$ we have (in the notation of [6])
$$
\Phi_0\mu_0F_0^{-3/2} \ll (nT)^{-3/4},
$$
which makes a total contribution of $O(1)$, as does error term
integral in Theorem 2.2 of [6]. The error terms with
$\Phi(a),\,\Phi(b)$ vanish for $a \to 0+,\,b \to +\infty\,$, and (3.1)
follows. Finally note that by using Taylor's formula
it is seen that the error made by replacing
$$
\left(\left({T\over2\pi n} + {1\over4}\right)^{1/2}
- {1\over2}\right)^{-1/2}
\exp\left(-G^2\left({\roman {arsinh}}\,\sqrt{\pi n\over 2T}\,\right)^{2}\right)
$$
with
$$
\left({T\over2\pi n}\right)^{-1/4}\exp\Bigl(-{\pi G^2n\over2T}\Bigr)
$$
in (3.1) is $\ll 1$ for $G \ge T^{1/5}\log^CT$.

\medskip
We remark that the series in (3.1) can be truncated at $T^{1/3}$. Namely the
contribution for $n \le T^{1/3}$ is, by trivial estimation,
$$
\ll T^{-1/4}\sum_{n\le T^{1/3}}d(n)n^{-1/4} \ll \log T,
$$
and this is absorbed by the $O$-term in (3.1).

%% \medskip
%% LEMMA 2. {\it If $A\in\RR$ is a constant, then we have}
%%$$
%%\cos\left(\sqrt{8\pi nT} +
%%\txt{1\over6}\sqrt{2\pi^3}n^{3/2}T^{-1/2} + A\right)
%%= \int_{-\infty}^\infty \a(u)\cos(\sqrt{8\pi n}(\sqrt{T} + u)
%%+ A)\d u,\leqno(3.4)
%% $$
%% {\it where $\a(u) \ll T^{1/6}$ for $u\not=0$,}
%% $$
%% \a(u) \ll T^{1/6}\exp(-bT^{1/4}|u|^{3/2}) \leqno(3.5)
%% $$
%% {\it for $u<0$, and}
%% $$
%% \a(u) = T^{1/8}u^{-1/4}\left(d\exp(ibT^{1/4}u^{3/2})
%% + {\bar d}\exp(-ibT^{1/4}u^{3/2})\right) + O(T^{-1/8}u^{-7/4})\leqno(3.6)
%% $$
%% {\it for $u \ge T^{-1/6}$ and some constants $b\; (>0)$ and $d$.}

%% \bigskip
%% This is a result of M. Jutila [21, Part II]. The formulas (3.5)--(3.6) enable
%% us to get rid of the term $n^{3/2}T^{-1/2}$ in $f(T,n)$ in Lemma 1. By using the
%% Taylor expansion, or a variation of the above method of Jutila, the influence of the
%% remaining terms in $f(T,n)$ can be made innocuous as well in our applications.

\medskip
LEMMA 2. {\it Let ${\Cal N}$ denote the number of solutions  in integers
m,n,k of the inequality}
$$
|\sqrt{m} + \sqrt{n} - \sqrt{k}\,| \;\le\; \delta\sqrt{M} \qquad(\delta > 0)
$$
{\it with $M' < n \le 2M', M < m \le 2M, k\in \NN$, and $M' \le M$. Then}
$$
{\Cal N} \;\ll_\e\; M^{\e}\Bigl(M^2M'\delta + (MM')^{1/2}\Bigr).\leqno(3.4)
$$

\medskip
LEMMA 3.  {\it Let $k\ge 2$ be a fixed
integer and $\delta > 0$ be given.
Then the number of integers $n_1,n_2,n_3,n_4$ such that
$N < n_1,n_2,n_3,n_4 \le 2N$ and}
$$
|n_1^{1/k} + n_2^{1/k} - n_3^{1/k} - n_4^{1/k}| \;<\; \delta N^{1/k}
$$
{\it is, for any given $\e>0$,}
$$
\ll_\e\; N^\e(N^4\delta + N^2).\leqno(3.5)
$$

\medskip
Lemma 2 was proved by Sargos and the author [18], while Lemma 3 is
due to Robert--Sargos [24]. They represent powerful arithmetic
tools which are essential in the analysis when the cube or
biquadrate of exponential sums involving $\sqrt n$ appears.

\medskip
LEMMA 4. {\it For $\,HU \gg T^{1+\e}\,$ and $\,T^\e \ll U \le \hf\sqrt{T}\,$
we have}
$$
\eqalign{& \int\limits_T^{T+H}\Bigl(E(x+U)-E(x)\Bigr)^2\d x \cr&=
{1\over4\pi^2}\sum_{n\le {T\over2U}}{d^2(n)\over
n^{3/2}}\int\limits_T^{T+H} x^{1/2}\left|\exp\left(2\pi iU\sqrt{{n\over
x}}\,\right)-1\right|^2\d x + O_\e(T^{1+\e} +
HU^{1/2}T^\e).\cr}\leqno(3.6)
$$

\medskip
This result was proved by M. Jutila [20]. The analogous formula   also holds
with $E(T)$ replaced by
$$
\D(x) \;:=\; \sum_{n\le x}d(n) -x(\log x + 2\gamma-1),\leqno(3.7)
$$
the error term in the classical Dirichlet divisor problem.
From  (3.6) Jutila deduced ($a\asymp b$ means that  $a\ll b\ll a$)
$$
\int_T^{2T}\Bigl(E(x+U)-E(x)\Bigr)^2\d x
\;\asymp\; TU\log^3\Bigl({\sqrt{T}\over U}\Bigr)
\qquad\Bigl(T^\e \ll U \le \hf\sqrt{T}\,\Bigr).
\leqno(3.8)
$$
The author   sharpened (3.8) to an asymptotic formula. Namely it
was proved in [15] that, with suitable constants $e_j\;(e_3>0)$ and
$T^\e \ll U \le \hf\sqrt{T}$,
$$\eqalign{
\int_T^{2T}\Bigl(E(x+U)-E(x)\Bigr)^2\d x & = TU\sum_{j=0}^3e_j\log^j
\Bigl({\sqrt{T}\over U}\Bigr) \cr&
+ O_\e(T^{1/2+\e}U^2) + O_\e(T^{1+\e}U^{1/2}).\cr}
$$

\head
4. The proof of Theorem 1
\endhead
We begin with the bound in (2.3). The left-hand side equals,
by the defining relation of $E(T)$ ((1.10)--(1.11) with $k=1$),
$$
\eqalign{&\int_T^{2T}|\zt|^2\Bigl(O(G\log T) + E(t+G)-E(t-G)\Bigr)\d t\cr&
\ll GT\log^2T + \int_T^{2T}|\zt|^2|E(t+G)-E(t-G)|\d t\cr&
\ll GT\log^2T +
{\Biggl(\int_T^{2T}|\zt|^4\d t \int_T^{2T}\Bigl(E(t+G)-E(t-G)\Bigr)^2
\d t\Biggr)}^{1/2}\cr&
\ll GT\log^2T + (T\log^4T\cdot TG\log^3T)^{1/2} \ll GT\log^2T\cr}
$$
for $G\ge T^\e$, as asserted. Here we used the Cauchy-Schwarz inequality
for integrals and (3.8). Note that the upper bound in (2.3) is best possible,
as it coincides with the lower bound in (1.5). An interesting, but difficult
problem, would be to obtain an asymptotic formula for the integral in (2.3).

\medskip
To discuss (2.4), we first exchange the order of integration in the
relevant integrals.
It follows that the left-hand side of (2.4) does not exceed
$$
\eqalign{&
\int_{T-G}^{2T+G}|\zx|^4\left(\int_{x-G}^{x+G}|\zt|^2\d t\right)\d x\cr&
= \int_{T-G}^{2T+G}|\zx|^4\Bigl(O(G\log T) + E(x+G) - E(x-G)\Bigr)\d x\cr&
\ll_\e\, GT\log^5T + T^{547/416+\e},\cr} \leqno(4.1)
$$
which immediately gives the bound which is somewhat weaker than the
one in (2.4), since $13/10 = 1.3 < 547/416 =  0.314903\ldots\,.$
Here we used the sharpest
known bound $E(T) \ll_\e T^{131/416+\e},\;131/416 = 0.314903\ldots\,$
of N. Watt [26]. To obtain the sharper bound asserted by (2.4)
we shall use results on the moments of $E^*(t)$ (see Section 5), and hence
the proof of the bound in question will be completed there.

\bigskip
For the proof of (2.5) we start from (1.2) which gives, for $T/2 \le t \le 5T/2$,
$$
|\zt|^2 \ll \log T\int_{t-T^\e}^{t+T^\e}|\zx|^2\d x + 1,
$$
and we use the trivial inequality
$$
\eqalign{&
\int_{t-G}^{t+G}|\zx|^k\d x = \int_{-G}^{G}|\z(\hf +it+iu)|^k\d u\cr&
\le {\roman e}\int_{-\infty}^\infty|\z(\hf +it+iu)|^k{\roman e}^{-(u/G)^2}\d u
= \sqrt{\pi}{\roman e}GJ_k(t,G)\cr}
$$
in the notation of (1.7), where $T/2 \le t \le 5T/2, 1 \ll G \ll T$
and $k\in\NN$ is fixed. This gives
$$
\eqalign{&
\int_T^{2T}|\zt|^2{\left(\int_{t-G}^{t+G}|\zx|^2\d x\right)}^2\,\d t\cr&
\ll T^\e\log T\int_{T/2}^{5T/2}\f(t)J_1(t,T^\e)
{\left(\int_{-\infty}^\infty|\z(\hf +it+iu)|^2{\roman e}^{-(u/G)^2}\d u
\right)}^2\d t\cr& \;+ G^2T\log^4T,\cr}
$$
following the proof of (2.3), where $\f(t)\,(\ge0)$ is a smooth function
supported in $[T/2,\, 5T/2]$, such that $\f(t) = 1$ for $T\le t\le2T$ and
$\f^{(r)}(t)\ll_r T^{-r}$ for $t\in\RR$ and any $r\in\NN$.
For $J_1(t,T^\e)$ we use Lemma 1,
writing $\sin z = ({\roman e}^{iz} - {\roman e}^{-iz} )/(2i)$, and integrate by
parts $\exp(i2\sqrt{2\pi nt}\,)$. In this way it is seen that
$$
\eqalign{&
\int_{T/2}^{5T/2}\f(t)\sum_{n=1}^\infty(-1)^nd(n)n^{-1/2}\ldots
\exp\bigl(if(t,n)\bigr){\left(\int_{-\infty}^\infty\ldots\right)}^2\d t\cr&
= \int_{T/2}^{5T/2}
\Biggl\{\f(t)\sqrt{t}\sum_{n=1}^\infty {i\over2\sqrt{\pi n}}(-1)^nd(n)n^{-1/2}\ldots\cr&
\times
\exp\Bigl( -i\txt{1\over4}\pi +
i\txt{1\over6}\sqrt{2\pi^3}n^{3/2}t^{-1/2} + a_5in^{5/2}t^{-3/2} +
\ldots\Bigr){\left(\int_{-\infty}^\infty\ldots\right)}^2\Biggr\}'\d t.
\cr}
$$
Note that, for $G = T^\e, t\asymp T$, we have
$$
\eqalign{&
{\left\{\exp\left(-G^2\left({\roman {arsinh}}
\sqrt{\pi n\over2t}\,\right)^2\right)\right\}}'\cr&
= {G^2\sqrt{\pi n\over2}\,{\roman {arsinh}}\sqrt{\pi n\over2t}\over t^{3/2}
\sqrt{1 + {\pi n\over2t}}}
\exp\left(-G^2\left({\roman {arsinh}}\sqrt{\pi n\over2t}\,\right)^2\right)\cr&
\asymp {G^2n\over t^2}\exp\left(-G^2\left({\roman {arsinh}}
\sqrt{\pi n\over2t}\,\right)^2\right),\cr}\leqno(4.2)
$$
$$
{(\f(t)\sqrt{t})}' \ll {1\over\sqrt{T}}\,,
$$
and that (here $G\ne T^\e$)
$$
\eqalign{&
\left\{{\left(\int_{-\infty}^\infty
|\z(\hf +it+iu)|^2{\roman e}^{-(u/G)^2}\d u \right)}^2\right\}'\cr&
= 2\int_{-\infty}^\infty|\z(\hf +it+iu)|^2{\roman e}^{-(u/G)^2}\d u
\int_{-\infty}^\infty{\d\over \d t}|\z(\hf +it+iu)|^2{\roman e}^{-(u/G)^2}\d u.
\cr}
$$
Integrating by parts we have
$$
\eqalign{&
\int_{-\infty}^\infty\Bigl\{{\d\over \d t}|\z(\hf +it+iu)|^2\Bigr\}
{\roman e}^{-(u/G)^2}\d u\cr&
= \int_{-\infty}^\infty\Bigl\{{\d\over \d u}|\z(\hf +it+iu)|^2\Bigr\}
{\roman e}^{-(u/G)^2}\d u\cr&
= 2\int_{-\infty}^\infty  uG^{-2}|\z(\hf +it+iu)|^2{\roman e}^{-(u/G)^2}\d u.\cr}
\leqno(4.3)
$$
Observe the integrals in (4.3) can be truncated at $|u| = G\log T$ with a negligible
error. Therefore, after an integration by parts, we get an integral with the
same type of exponential factor (i.e., $f(t,n)$ in the exponential), but there
will be in the integrand a smooth factor of the order
$\ll G^{-1}\sqrt{T/n}$. Hence after a large number of integrations by parts it follows
that the contribution of $n$ satisfying $n > T^{1+\e}G^{-2}$ will be negligible (i.e.,
less than $T^{-A}$ for any given $A>0$ and $\e = \e(A)$).
This truncation of the series over $n$
is the crucial point in the proof, as the ensuing expression will be quite
similar to the expressions for $J_1(t,G)$, only in the
exponential factor in (4.2) we shall have $G = T^\e$. Thus the proof reduces
to the estimation of
$$
T^\e G^2\int_{T/2}^{5T/2}\f(t)\sum\nolimits_1{\left(\sum\nolimits_2\right)}^2\d t,
\leqno(4.4)
$$
where
$$\eqalign{
\sum\nolimits_2 &:=
\sum_{n\le T^{1+\e}G^{-2}} (-1)^nd(n)n^{-1/2}
\left({\Bigl({t\over2\pi n} + {1\over4}\Bigr)}^{1/2}
- {1\over2}\right)^{-1/2}\times\cr&\times\exp\left(-G^2{\Bigl(
{\roman {arsinh}}\,\sqrt{\pi n/(2t)}\,\Bigr)}^{2}\right)\sin f(t,n),
\cr}
$$
and $\sum\nolimits_1$ is the same expression with $G=T^\e$ in the
exponential factor. The other two terms, which arise after the
squaring of the right-hand side of (3.1), are clearly less difficult
to deal with. Note that
$$\eqalign{&
\sum\nolimits_1{\left(\sum\nolimits_2\right)}^2
= \sum_{m\le T^{1+\e}G^{-2}}(-1)^md(m)\cdots\sin f(t,m)\times\cr&
\sum_{n\le T^{1+\e}G^{-2}}(-1)^nd(n)\cdots \sin f(t,n)
\sum_{k\le T^{1+\e}G^{-2}}(-1)^kd(k)\cdots \sin f(t,k),\cr}
$$
and write the sines as exponentials. For $G\ge T^{1/7+\e}$ we use Taylor's formula
to remove the terms $a_7n^{7/2}T^{-5/2}+\ldots$ from all functions $f$ coming from Lemma 1.
Namely we can truncate the tails of series after sufficiently many terms to obtain a negligible
error term. There remain only finitely many terms, but the exponentials are identical, so
it suffices to treat the first terms only.
Then we integrate by parts many times, as was done in the previous
part of the proof. Thus we are left with sums containing the exponential
$$
e^{i(\D t^{1/2} + Et^{-1/2} + Ft^{-3/2})} = e^{if(t)},
$$
say, where we set
$$\eqalign{
\D &:= \sqrt{8\pi}(\sqrt{m}+\sqrt{n}-\sqrt{k})
\cr
E &:= \txt{\frac{1}{6}}\sqrt{2\pi^3}(m\sqrt{m}+n\sqrt{n}-k\sqrt{k})
\cr
F &:= a_5(m^2\sqrt{m}+n^2\sqrt{n}-k^2\sqrt{k}).
\cr}\leqno(4.5)
$$
Namely the terms with
$$
\sqrt{m} + \sqrt{n} + \sqrt{k}, \quad-\sqrt{m} - \sqrt{n} - \sqrt{k}
$$
are clearly negligible by sufficiently many integrations by parts.
Thus only the combination of signs as in (4.5)
is relevant.
Here we suppose that
$$
M' < n \le 2M', \quad M < m \le 2M, \quad K \le k \le 2K, \quad  M'\le M,\quad M,N,K \ge T^{1/3},
$$
and consider first the contribution
 from the triplets $(m,n,k)\;(\in\NN^3)$ satisfying $\D \le T^{\e-1/2}$.
 We suppose $\D>0$, since the case $\D<0$ is analogous, and the case $\D=0$ is easy.
  Furthermore, by the first derivative test ([6, Lemma 2.1])
it is seen that the contribution is small if $K < AM$ or $K > BM$
with suitable positive constants $A,B$ (when $\D \gg \sqrt{M}$ or $\D \gg \sqrt{K}$).
Therefore, by using the bound  (3.4) of Lemma 2
(with $\delta = T^{\e-1/2}M^{-1/2}$) it is seen that
the corresponding portion of the integral in (4.4) is
$$
\eqalign{&
\ll_\e T^{1+\e}\max_{K,M,M'\le T^{1+\e}G^{-2},K\asymp M}
T^{-3/4}M^{-3/4}(M^{5/2}T^{-1/2} + M) + T^{1+\e}\cr&
\ll_\e T^{1+\e}\max_{K,M,M'\le T^{1+\e}G^{-2},K\asymp M}
(T^{-5/4}M^{7/4} + T^{-3/4}M^{1/4})+ T^{1+\e}\cr&
\ll_\e T^{1+\e}(T^{1/2}G^{-7/2} + 1) \ll_\e T^{1+\e}\cr}
$$
for $G\ge T^{1/7+\e}$, as asserted.

\medskip
Now we proceed analogously as was done in the author's work [14]. Suppose
$\D \ge T^{\e-1/2}$. We may  assume that $E >0$, since the other case
is analogous. Let
$$
T^{\e-1/2} \;\leqslant\;\D\;\leqslant\; \D_0,
$$
where $\D_0$, which will be determined later, does not depend on $m,n,k$.
Further suppose  that
$$2^{-j}\D_0 < \D \leqslant 2^{1-j}\D_0\qquad\Bigl(1 \leqslant j\leqslant J
\;(\asymp \log(\D_0T^{1/2-\e})\Bigr).
$$
If $|F|T^{-3/2} \ll ET^{-1/2}$ or $|F|T^{-3/2} \ll \D T^{1/2}$ with suitable
$\ll$-constants, then in $e^{if(t)}$ either $\D t^{1/2}$ or $E t^{-1/2}$
dominates in size. Hence we can use the method of [14].
If we have $\D > C_1E/T$ with
a sufficiently large $C_1>0$, then
$f'(t) \gg \D/\sqrt{T}$ in $[T,3T]$.
 Also
if $\D \geqslant T^{\e-1/2}$ and $\D < C_2E/T$ with a sufficiently small $C_2>0$,
then $f'(t) \gg ET^{-3/2}$. In both cases we estimate the integral of
$e^{if(t)}$ by the first derivative test, and then the sum over $m,n,k$ by Lemma 3.

\medskip
If $\D \geqslant T^{\e-1/2}$ and $\D \asymp E/T$, then there
may exist a saddle point $t_0 = E/\D$ (root of $f'(t_0) =0$)
in $[T,\,3T]$ if $\D T \asymp E$. Hence by the saddle-point method
(see [6, Chapter 2] or by the use the second derivative test,
making first the change of
variable $\sqrt{t} = u$, we obtain
$f''(t_0) = \hf \D^{5/2}E^{-3/2} \asymp \D T^{-3/2}$. Hence by the second
derivative test (see Lemma 2.2 of [6])
the corresponding portion of the integral in (4.4) is
$$
\eqalign{&
\ll
\sum_{j\le J}2^{j/2}\D_0^{-1/2}T^{3/4}T^{-3/4}M^{-3/4}(M^{5/2}\D_02^{-j} +M)
\cr&
\ll T^{\e}(M^{7/4}\D_0^{1/2} + M^{1/4}T^{1/4} +T) \ll_\e T^{1+\e}
\cr}
$$
for $M^{7/4}\D_0^{1/2}\ll T$, or
$$
\D_0 \ll T^2M^{-7/2}.
$$
But since trivially $\D_0 \ll \sqrt{M}$,  and $M\ge T^{1/3}$ it follows that
$$
T^2M^{-7/2} \ge T^{5/6} > \sqrt{M} \gg \D_0,
$$
which is needed.

\medskip
Finally if $|F|T^{-3/2} \gg ET^{-1/2}$ and $|F|T^{-3/2} \gg \D T^{1/2}$, then
$$
f'(t) \gg |F|T^{-5/2},\quad \D \ll |F|T^{-2}.
$$
The first derivative test shows that the contribution is, since $M \le T^{1+\e}G^{-2}$,
$$
\eqalign{&
\ll
T^{5/2}|F|^{-1}M^{-3/4}T^{-3/4}(M^{5/2}|F|T^{-2} +M)
\cr&
\ll T^{-1/4}M^{7/4} + M^{1/4}T^{7/4}|F|^{-1}
\cr&
 \ll_\e T^{1+\e}(T^{1/2}G^{-7/2}+ M^{1/4}T^{3/4}|F|^{-1}) \ll_\e T^{1+\e}
\cr}
$$
for $G\ge T^{1/7+\e}$, provided that
$$
M^{1/4}T^{3/4}|F|^{-1} \ll 1.
$$
If, however, $|F| \ll M^{1/4}T^{3/4}$, then $\D T^2 \ll |F| \ll M^{1/4}T^{3/4}$ and this implies
$$
\D \ll M^{1/4}T^{-5/4} <  T^{-1/2},
$$
and this  case has been already dealt with. This completes the proof of (2.5).

\bigskip
The proof of (2.6) is similar to the proof of (2.5).
The major difference is that, instead of (4.4), now we shall have to bound
$$
T^\e G^3\int_{T/2}^{5T/2}\f(t)\sum\nolimits_1{\left(\sum\nolimits_2\right)}^3\d t.
\leqno(4.6)
$$
We use then H\"older's inequality to deduce that the integral in (4.6)
does not exceed
$$
\left(\int_{T/2}^{5T/2}\f(t)\Bigl|\sum\nolimits_1\Bigr|^4\d t\right)^{1/4}
\left(\int_{T/2}^{5T/2}\f(t)\Bigl|\sum\nolimits_2\Bigr|^4\d t\right)^{3/4}.
\leqno(4.7)
$$
Both integrals in (4.7) are estimated similarly. Here we have
$$
T^{{1\over5}+\e} \ll G = G(T) \ll T^{1/3}.
$$
Therefore, by using Taylor's theorem, instead of the exponential
$$
e^{i(\D t^{1/2} + Et^{-1/2} + Ft^{-3/2})}
$$
we shall have the simpler function $e^{i(\D t^{1/2} + Et^{-1/2})}= e^{iH(t)}$, say.

\medskip
First note that the sum over $n$ in $\sum_1$ is split
into $O(\log T)$ subsums where $M < n\le M'\le 2M$, with $M \ll T^{1+\e}G^{-2}$.
Instead of Lemma 2 we use (3.5) of Lemma 3
(with $\delta = \D M^{-1/2}$),
supposing first that $\D>0$ and that
 $\D \le T^{\e-1/2}$. Afterwards the integral is estimated trivially.
The contribution to the relevant integral in (4.7) will be
$$\eqalign{&
\ll_\e T^{1+\e}\max_{M\ll T^{1+\e}G^{-2}}T^{-1}M^{-1}(M^4T^{\e-1/2}M^{-1/2}
+ M^2) + T^{1+\e}\cr& \ll_\e T^{1+\e}(TG^{-5}+1) \ll_\e T^{1+\e}\cr}
$$
for $G \ge T^{1/5+\e}$, as asserted.

\medskip
Now suppose $\D > T^{\e-1/2}$.
If $\D>0$, we may also assume that $E>0$, for otherwise all derivatives
of $H(t)$ have the same sign. Let $\D > C_1E/T$ for some suitable
$C_1>0$. Then $H'(t) \gg \D T^{-1/2}$, and supposing that
$$
\D \asymp 2^jT^{\e-1/2} \qquad(j = 1,2,\ldots),
$$
we obtain by the first derivative test that the contribution is
$$
\eqalign{&
\ll
\sum_{j\ge 1}M^{-1}T^{-1}2^{-j}T^{-\e}(M^{7/2}2^jT^{\e-1/2} + M^2)
\cr&
\ll_\e M^{5/2}T^{\e-3/2} + T^\e \ll_\e T^{1+\e}.
\cr}
$$
Similar arguing is used if $\D < C_2E/T$. There remains the case when $\D \asymp E/T$,
in which case (after the substitution $u = \sqrt{t}$) it  is seen that
$$
|H(t)|^{-1/2} \asymp \D T^{-3/2}.
$$
Let, as in the proof of (2.5),
$$2^{-j}\D_0 < \D \leqslant 2^{1-j}\D_0\qquad\Bigl(1 \leqslant j\leqslant J
\;(\asymp \log(\D_0T^{1/2-\e})\Bigr).
$$
For $M \le T^{1/2}$ the contribution is
$$
\eqalign{&
\ll
\sum_{j\le J}M^{-1}T^{-1}2^{j/2}\D_0^{-1/2}T^{3/4}(M^{7/2}2^{-j}\D_0 + M^2)
\cr&
\ll_\e M^{5/2}T^{-1/4}\D_0^{1/2} + MT^\e \ll_\e T^{1+\e},
\cr}
$$
for $\D_0 \le T^\e$, when
$$
M^{5/2}T^{-1/4}\D_0^{1/2} \le T\D_0^{1/2} \ll_\e T^{1+\e/2}.
$$
But if $\D_0\ge T^\e$, then since $\D \asymp E/T$, successive integrations
by part of $e^{i\D\sqrt{t}}$ show that the contribution is negligible.

\medskip
Let now
$$
T^{1/2} < M \le T^{1+\e}G^{-2} \le T^{3/5}\qquad(G\ge T^{1/5+\e}).
$$
Like in the preceding case, the contribution will be
$$
\ll_\e M^{5/2}T^{-1/4}\D_0^{1/2} + MT^\e \ll_\e T^{1+\e}
$$
for
$$
\D_0 \le T^{5/2+\e}M^{-5}.\leqno(4.8)
$$
If (4.8) does not hold, then we have $\D_0 > T^\e$, since $M > T^{1/2}$.
Again we integrate by parts $e^{i\D\sqrt{t}}$ suffieciently many times.
Each time we get a factor
in the integrand which is
$$
\eqalign{&
\ll \frac{1}{\sqrt{T}\D} + \frac{\sqrt{T}}{\D} \cdot T^{-3/2}E
\cr&
\ll T^{-1/2-\e} + T^{-\e-1}M^{3/2}
\cr&
\ll T^{-1/2-\e} + T^{-\e-1}T^{9/10} < T^{-1/10},\cr}
$$
so that the contribution will be negligible. This proves (2.6).

\medskip
To prove (2.7), note that the left-hand side equals (see (1.11) with $k=2$)
$$\eqalign{&
\int_T^{2T}|\zt|^2\Bigl(O(G\log^4T) + E_2(t+G)-E_2(t-G)\Bigr)^2\d t\cr&
\ll TG^2\log^{9}T + {\Bigl(\int_T^{2T}|\zt|^4\d t
\int_{T/3}^{3T}E_2^4(t)\d t\Bigr)}^{1/2}\cr&
\ll TG^2\log^{9}T + T^{13/6}\log^CT \ll TG^2\log^{12}T
\cr}\leqno(4.9)
$$
for $G \ge T^{7/12}\log^CT$, as asserted. Here we used the bound, which
follows from (1.14), namely
$$
\int_0^T|E_2(t)|^A\d t \ll T^{2+{2\over3}(A-2)}\log^CT
\quad(A\ge 2,\,C = 22 + 8(A-2))
\leqno(4.10)
$$
with $A=4$. This completes the proof of Theorem 1.

\bigskip
\head
5. The proof of Theorem 2
\endhead
We have first, similarly as in (4.9),
$$
\eqalign{&
\int_T^{2T}|\zt|^4\int_{t-G}^{t+G}|\zx|^4\d x\,\d t\cr&
= \int_T^{2T}|\zt|^4\Bigl(O(G\log^4 T) + E_2(t+G) - E_2(t-G)\Bigr)\d t\cr&
\ll TG\log^8T + T^{5/3}\log^CT,\cr}\leqno(5.1)
$$
where we used  the first bound in (1.14).
This implies that the left-hand side of (5.1) is $\ll (TG+T^{5/3})\log^CT$ in
the whole range $1 \ll G = G(T) \ll T$. The bound in question was actually proved
by Ivi\'c-Jutila-Motohashi [19] in connection with mean square estimates
for ${\Cal Z}_2(s)$. Hence the main task is to prove the other bound
in (2.8), for which we need the second bound in (1.14). Note
that the left-hand side of (2.8) is majorized by a multiple ($L = \log T$) of
$$
\eqalign{&
\int_T^{2T}|\zt|^4\int_{-\infty}^\infty|\z(\hf+it+iu)|^4{\roman e}^{-(u/G)^2}
\d u \d t\cr&
\ll TGL^8 + \int_T^{2T}|\zt|^4\int_{-\infty}^\infty
\left\{{\d\over\d t}E_2(t+u)\right\}
{\roman e}^{-(u/G)^2}\d u \d t\cr&
\ll TGL^8 + 2\int_T^{2T}|\zt|^4\int_{-\infty}^\infty E_2(t+u)uG^{-2}
{\roman e}^{-(u/G)^2}\d u \d t\cr&
\ll TGL^8 + {L\over G}\int_{T-GL}^{2T+GL}|E_2(u)|\int_{u-GL}^{u+GL}|\zt|^4\d t\d u\cr&
\ll L^C\left\{TG + T^{3/2} + G^{-2}\int_{T-GL}^{2T+GL}|E_2(u)|
\int_{u-GL^2}^{u+GL^2}|E_2(t)|\d t\d u\right\}.\cr}\leqno(5.2)
$$
Here we used the fact that
$$
{\d\over\d t}E_2(t+u) = {\d\over\d u}E_2(t+u),\leqno(5.3)
$$
and integrated by parts. We used a similar procedure later, namely
($Q_4 = P_4 + P'_4$ (see (1.11)--(1.12))
is a suitable polynomial of degree four, $C$ is henceforth
a generic positive constant, $u\asymp T$)
$$
\eqalign{&
\int_{u-GL}^{u+GL}|\zt|^4\d t \ll
\int_{-\infty}^\infty|\z(\hf+iu+iv)|^4{\roman e}^{-(v/GL)^2}\d v
\cr&
= \int_{-\infty}^\infty \Bigl(Q_4(\log(u+v)) + E'_2(u+v)\Bigr)
{\roman e}^{-(v/GL)^2}\d v\cr&
= O(GL^C) + 2\int_{-\infty}^\infty vG^{-2}E_2(u+v){\roman e}^{-(v/GL)^2}\d v
\cr&\ll L^C\Bigl(G + G^{-1}\int_{u-GL^2}^{u+GL^2}|E_2(t)|\d t\Bigr).\cr}
$$
We also used the bound, which follows by the Cauchy-Schwarz inequality
for integrals from the mean square bound in (1.14), namely
$$
\int_0^T|E_2(t)|\d t \;\ll\;T^{3/2}\log^CT.
$$
To complete the proof it remains to note that
$$
\eqalign{&
{\left(\int_{T-GL}^{2T+GL}|E_2(u)|
\int_{u-GL^2}^{u+GL^2}|E_2(t)|\d t\d u\right)}^2\cr&
\ll \int_{T-GL}^{2T+GL}E_2^2(u)\d u \int_{T-GL}^{2T+GL}
{\Bigl(\int_{u-GL^2}^{u+GL^2}|E_2(t)|\d t\Bigr)}^2\d u\cr&
\ll T^2L^C \int_{T-GL}^{2T+GL} GL^2 \int_{u-GL^2}^{u+GL^2}E_2^2(t)\d t\d u\cr&
\ll T^2GL^C \int_{T-GT^\e}^{2T+T^\e}E_2^2(t)\left(\int_{t-GL^2}^{t+GL^2}\d u
\right)\d t\cr&
\ll T^4L^CG^2.\cr}\leqno(5.4)
$$
When we take the square root in (5.4) and insert the resulting bound in (5.2)
we are left with the bound
$$
O\left(L^C(TG+ T^{3/2} + T^2G^{-1})\right)
$$
for the left-hand side of (2.8).  But as $T^{3/2} \le TG$ for $G\ge T^{1/2}$
and $T^{3/2} \le T^2G^{-1}$ for $G \le T^{1/2}$, this means that in
the bound above the term
$T^{3/2}$ may be omitted, and (2.8) follows.
We point out yet another estimate, namely (6.2), for the integral in (2.8).
This was derived for the proof of Theorem 3, and does not contain the
(expected) term $TG\log^CT$, but terms which are reasonably small when
$G$ is `about' $\;T^{131/416}$.

\bigskip
The proof of (2.9) is based on the use of (5.3) and the fourth moment of the
function $E^*(t)$, defined by
$$
E^*(t) \;:=\; E(t) - 2\pi\D^*\bigl({t\over2\pi}\bigr),
$$
where (see (3.7))
$$
\D^*(x) := -\D(x)  + 2\D(2x) - \hf\D(4x)
= \hf\sum_{n\le4x}(-1)^nd(n) - x(\log x + 2\gamma - 1),
$$
and $\D(x)$ is the error term in the  Dirichlet divisor problem.
The function $E^*(t)$ was investigated by several authors, including
M. Jutila [21], who introduced it, and the author [10], [11], [12] and [14].
Among other things, the author ([10, Part II] and [14]) proved that
$$
\int_0^T|E^*(t)|^5 \d t \;\ll_\e\; T^{2+\e}\leqno(5.5)
$$
and that ([10, Part IV, Corollary 2] and [14])
$$
\int_0^T\bigl(E^*(t)\bigr)^4 \d t \;\ll_\e\; T^{7/4+\e}.\leqno(5.6)
$$
The advantage of working with $E^*(t)$ instead of $E(t)$ is that the
former is, in the mean power sense, smaller than the latter (for this
see [6, Chapter 15] and [10]).

\medskip
For our proof  we need from [13] the elementary formula (cf. (3.4))
$$
J_1(t,G) = {2\over\sqrt{\pi}G^3}\int_{-G\log T}^{G\log T}
xE^*(t+x){\roman e}^{-(x/G)^2}\d x + O(\log^2 T),\leqno(5.7)
$$
which is valid for $T^\e \le G = G(T) \le T^{1/3},\,T/2\le t \le 5T/2$.
First observe that the left-hand side of (2.9) is majorized ($\f(t)$ is as in
the proof of (2.5)) by
$$
\eqalign{&
G^2\int_{T/2}^{5T/2}|\zt|^4\f(t)J_1^2(t,G)\d t\cr& =
G^2\int_{T/2}^{5T/2}(Q_4(\log t) + E_2'(t))\f(t)J_1^2(t,G)\d t\cr&
= O(TG^2\log^8T) - G^2\int_{T/2}^{5T/2}E_2(t)
{\Bigl(\f(t)J_1^2(t,G)\Bigr)}'\d t.
\cr}\leqno(5.8)
$$
Namely by the Cauchy-Schwarz inequality for integrals and the classical bound
$$
\int_0^T|\zt|^4\d t \ll T\log^4T,
$$
we have
$$\eqalign{&
G^2\int_{T/2}^{5T/2}Q_4(\log t)\f(t)J_1^2(t,G)\d t \cr&
\ll G\log^4T\int_{T/2}^{5T/2}\int_{-\infty}^\infty|\z(\hf+it+iu)|^4
{\roman e}^{-(u/G)^2}\d u \d t\cr&
\ll G\log^4T\int_{T/3}^{3T}|\zx|^4\int_{x-GL}^{x+GL}
{\roman e}^{-(x-t)^2/G^2}\d t\d x \ll TG^2\log^8T.\cr}
$$
We also use $\f'(t) \ll 1/T$ (the contribution of this derivative is easily
handled) and (5.7). Hence with suitable $g(t)$ we obtain
$$
{\Bigl(J_1^2(t,G)\Bigr)}' = 2J_1(t,G)\Bigl\{-{2\over\sqrt{\pi}G^3}
\int\limits_{-G\log T}^{G\log T}
E^*(t+x)\bigl(1 - {2x^2\over G^{2}}\bigr){\roman e}^{-(x/G)^2}\d x + g'(t)\Bigr\},
\leqno(5.9)
$$
where $g(t) \ll \log^2T$. The term $g'(t)$ is integrated back by parts,
and its contribution is easily seen to be $\ll_\e T^{1+\e}G^2$ in the
required range. The first term on the right-hand side of (5.9)
is inserted in (5.8). It follows that the main contribution to the left-hand
side of (5.8) is ($L = \log T)$ bounded by $T^{1+\e}G^2$ plus
$$
\eqalign{&
G^2\int_{T/2}^{5T/2}|E_2(t)|\f(t)G^{-5}\Bigl(\int_{-GL}^{GL}|E^*(t+x)|\d x\Bigr)^2\d t\cr&
\ll G^{-3}{\Bigl(\int_{T/2}^{5T/2}|E_2(t)|^2\d t\Bigr)}^{1/2}
{\Bigl(\int_{T/2}^{5T/2}\Bigl(\int_{-GL}^{GL}|E^*(t+x)|\d x\Bigr)^4\d t\Bigr)}^{1/2}\cr&
\ll TL^CG^{-3}
{\Biggl(\int_{T/2}^{5T/2}G^3\int_{t-GL}^{t+GL}|E^*(u)|^4\d u\d t\Biggr)}^{1/2}\cr&
\ll TL^CG^{-3}{\Biggl(\int_{T/2-GL}^{5T/2+GL}G^3|E^*(u)|^4
\Bigl(\int_{u-GL}^{u+GL}\d t\Bigr)\d u\Biggr)}^{1/2}
\cr&
\ll_\e T^{15/8+\e}G^{-1},\cr}
$$
where we used H\"older's inequality for integrals, (5.3) and (5.6).
Since
$$
T^{15/8}G^{-1} \le TG^2\qquad ({\roman {for}}\;
G \ge T^{7/24}),
$$
 the bound in (2.9) follows.

\bigskip
Finally it remains to prove the bound (2.10) of Theorem 2. We proceed similarly
as in the proof of (2.9), and we majorize  first $\int_{T}^{2T}\ldots$  by
$\int_{T/2}^{5T/2}\f(t)\ldots$. Then we majorize
$$
{\left(\int_{t-G}^{t+G}|\zx|^4\d x\right)}^2
$$
by $G^2J_2^2(t,G)$ and write
$$
|\zt|^4 = Q_4(\log t) + E_2'(t), |\z(\hf + it + iu)|^4 = Q_4(\log (t+u)) + E_2'(t+u).
$$
After this we integrate by parts $E_2'$, using (5.3) and
obtaining $E_2(t), E_2(t+u)$, but gaining
essentially a factor of $1/G$ each time in the process. The major contribution to the
left-hand side of (2.10) will be $\ll_\e T^{1+\e}G^2$ plus
$$\eqalign{&
G^{-3}\int_{T/2}^{5T/2}\f(t)|E_2(t)|{\left(\int_{-GL}^{GL}
|E_2(t+x)|^4{\roman e}^{-(x/G)^2}\d x\right)}^2\d t\cr&
\ll G^{-3}{\Bigl(\int_{T/2}^{5T/2}\f(t)|E_2(t)|^2\d t\Bigr)}^{1/2}
{\Biggl(\int_{T/2}^{5T/2}{\Bigr(\int_{t-GL}^{t+GL}
|E_2(u)|\d u\Bigr)}^4\d t\Biggr)}^{1/2}.
\cr}\leqno(5.10)
$$
Now we use (4.10) with $A=4$ and H\"older's inequality, to obtain that
$$\eqalign{&
\int_{T/2}^{5T/2}{\left(\int_{t-GL}^{t+GL}|E_2(u)|\d u\right)}^4\d t\cr&
\ll\ (GL)^3\int_{T/2}^{5T/2}\int_{t-GL}^{t+GL}|E_2(u)|^4\d u\cr&
= \;GL^3\int_{T/2-GL}^{5T/2+GL}E_2^4(u)\left(\int_{u-GL}^{u+GL}\d t\right)
 \d u\cr&
\ll \;T^{10/3}G^4\log^CT.\cr}\leqno(5.11)
$$
Hence if we insert (5.11) in (5.10) and use (5.3), we obtain that the expression
in (5.10) is
$$
\ll\; T^{8/3}G^{-1}\log^CT \;\ll_\e\; T^{1+\e}G^2
$$
for $G \ge T^{5/9}$, which yields then (2.10) and completes the proof of Theorem 2.

\medskip
It remains yet to complete the proof of (2.4). From (4.1) we have
$$
\eqalign{&
\int_{T-G}^{2T+G}|\zx|^4\left(\int_{x-G}^{x+G}|\zt|^2\d t\right)\d x\cr&
= \int_{T-G}^{2T+G}|\zx|^4\Bigl(O(G\log T) + E(x+G) - E(x-G)\Bigr)\d x\cr&
\ll_\e GT^\e + \int_{T-G}^{2T+G}|\zx|^4\Bigl|E^*(x+G) - E^*(x-G)\Bigr|\d x.
\cr} \leqno(5.12)
$$
Here we used the defining property of $E^*(t)$ together with the
elementary bound
$$
\D^*(x+G) - \D^*(x-G) \;\ll_\e\;GT^\e\qquad(x\asymp T,\;1\ll G\ll T),
$$
which is easily obtained, since
$$
\D^*(x) = \hf\sum_{n\le4x}(-1)^nd(n),\qquad d(n) \ll_\e n^\e.
$$
At this point we use H\"older's inequality for integrals, (5.5) and
the bound (see [6, Chapter 8])
$$
\int_0^T|\zt|^5\d t \;\ll_\e\;T^{9/8+\e},
$$
to deduce that the last integral in (5.12) is
$$
\eqalign{&
\le {\left(\int_{T/3}^{3T}|\zx|^5\d x\right)}^{4/5}
{\left(\int_{T/3}^{3T}|E^*(x)|^5\d x\right)}^{1/5}\cr&
\ll_\e T^{{9\over8}\cdot{4\over5}+2\cdot{1\over5}+\e} = T^{{13\over10}+\e},
\cr}
$$
which yields (2.4). For small values of $G$,
namely for $1\ll G \ll T^{1/10}$,
the bound in (2.4) may be improved by using a more general result than (3.8),
namely
$$
\int_T^{T+H}\Bigl(E(x+U)-E(x)\Bigr)^2\d x
\;\asymp\; HU\log^3\left({\sqrt{T}\over U}\right),
\leqno(5.13)
$$
deduced by M. Jutila [20] from (3.6), for $HU \gg T^{1+\e}$ and
$T^\e \ll U \le \hf\sqrt{T}$. If (5.13) is applied in
conjunction with (4.1), the Cauchy-Schwarz inequality
and the bound (see (2.11)) $\int_0^T|\zt|^8\d t \ll_\e T^{3/2+\e}$,
we obtain (2.4) with $T^\e(GT + T^{5/4}G^{1/2})$, and
$$
T^{5/4}G^{1/2} \;\le\; T^{13/10}\qquad(1\ll G \le T^{1/10}).
$$
This is unconditional, but conjectures on the order
of $E(T)$  and $E^*(T)$ would lead to further improvements,
e.g., the conjectural bound
$E(T) \ll_\e T^{1/4+\e}$ would replace the exponent 13/10 in (2.4) by 5/4.

\medskip
We end this section by pointing out that one can improve (1.6) for
the range given in Section 1. Namely, since the integral $J_k(t,G)$
in (1.7) can be truncated at $u = \pm G\log T$ with an error which is
$\ll T^{-A}$ for any given $A>0$, it follows that (1.6) is equivalent to
$$
\int_T^{2T}{\Bigl(\int_{t-G}^{t+G}|\zx|^4\d x\Bigr)}^m\d t
\ll_\e T^{1+\e}G^m.\leqno(5.14)
$$
We proceed as in the proof of (2.7), using (4.10), to infer that the
integral in (5.14) is (as before $C>0$ is a generic constant)
$$
\eqalign{&
\ll TG^m(\log T)^{4m} + \int_{T/3}^{3T}|E_2(t)|^m\d t\cr&
\ll TG^m(\log T)^{4m} + T^{2+{2\over3}(m-2)}\log^CT\cr&
\ll TG^m(\log T)^{C}
\cr}
$$
for
$$
G\;\ge\; T^{2m-1\over3m}\qquad(m\ge2),\leqno(5.15)
$$
where incidentally $m$ does not have to be an integer. In particular, it
follows then from (5.15) that (1.6) holds
for $G\ge T^{1/2}\;(m=2)$, $G \ge T^{5/9}
\;(m=3)$ and $G\ge T^{7/12}\;(m=4)$. Since $5/9 < 4/7$ and
$7/12 < 3/5$, this means that we have improved the range of $G$ for which
(1.6) holds when $m=3,4$.

\head
6. The proof of Theorem 3
\endhead
\no
The method of obtaining mean square bounds for $\Z_2(s)$ (see (1.9)) was
developed in [12] and [19], so that we shall be fairly brief.
From [19, pp. 337-339], we have
$$\eqalign{&
\int_T^{2T}|{\Cal Z}_2(\s + it)|^2\d t
\ll_\e
T^{2+\e}X^{1-2\s} +\cr& \,+ \log T\sup_{T^{1-\e}\le K\le X}
\int\limits_{T/2}^{5T/2}
\int\limits_K^{2K}|\zx|^4x^{-2\s}
\int\limits_{x-KT^{\e-1}}^{x+KT^{\e-1}}|\zy|^4\d y\d x\d t
\cr& \ll_\e
T^{2+\e}X^{1-2\s} + T\log T\sup_{T^{1-\e}\le K\le X}K^{-2\s}
\int\limits_K^{2K}|\zx|^4
\int\limits_{x-G}^{x+G}|\zy|^4\d y\d x\cr}\leqno(6.1)
$$
with $G = KT^{\e-1},\,\s> 1/2$, and $X$ a parameter to be suitably chosen.
For the last integral above one
could use (2.8) of Theorem 2. However, this would not lead to the result
of Theorem 3, as we need a bound when $G$  is `about' $K^{1/3}$,
or even slightly smaller. Thus
we shall proceed differently, and
use the elementary inequality
$$
ab \le \hf(a^{1/2}b^{3/2} + a^{3/2}b^{1/2})\qquad(a,b \ge 0)
$$
to obtain
$$
\eqalign{&
\int_{T}^{2T}|\zt|^4\int_{t-G}^{t+G}|\zx|^4\d x\d t
\cr&
\ll \int_{T}^{2T}\int_{t-G}^{t+G}\Bigl(|\zt|^2|\zx|^6+|\zt|^6|\zx|^2\Bigr)
\d x\d t\cr&
\ll \int_{T-G}^{2T+G}|\zx|^6\int_{x-G}^{x+G}|\zt|^2\d t \d x\cr&
+ \int_{T}^{2T}|\zt|^6\int_{t-G}^{t+G}|\zx|^2\d x\d t\cr&
\ll \int_{T/3}^{3T}|\zt|^6\Bigl(O(G\log T) + E(t+G)-E(t-G)\Bigr)\d t \cr&
\ll_\e
T^{\e+{1+\rho\over2}}(G + T^{131\over416}).
\cr}\leqno(6.2)
$$
Here we used the defining relation of $E(T)$ together
with the sharpest known bound $E(T) \ll_\e T^{131/416+\e}$ (see [5], [6]).
We also used the bound for the sixth
moment which follows from (2.12) and the Cauchy-Schwarz inequality,
namely
$$
\int_0^T|\zt|^6\d t \le {\left(\int_0^T|\zt|^4\d t\int_0^T|\zt|^8\d t
\right)}^{1/2} \ll_\e T^{{1+\rho\over2}+\e}.
$$
We use (6.2) (with $G = KT^{\e-1}$) in (6.1) to obtain

$$
\eqalign{&
T\sup_{T^{1-\e}\le K\le X}K^{-2\s}\int_K^{2K}|\zx|^4
\int_{x-G}^{x+G}|\zy|^4\d y\d x\cr&
\ll_\e T^{1+\e}\sup_{T^{1-\e}\le K\le X}\left(
K^{{1+\rho\over2}+{131\over416}-2\s} + T^{-1}K^{{3+\rho\over2}-2\s}\right)
\cr&
\ll_\e T^{1+\e} + T^{\e}X^{{3+\rho\over2}-2\s}
\cr}
$$
if
$$
 {1+\rho\over4} + {131\over832} \le \s \le {3+\rho\over4} .\leqno(6.3)
$$
If (6.3) holds, then  we obtain from (6.1)
$$\eqalign{
\int_T^{2T}|{\Cal Z}_2(\s + it)|^2\d t &\ll_\e T^{\e}(T
 + X^{{3+\rho\over2}-2\s} + T^2X^{1-2\s})\cr&
 \ll_\e
 T^\e(T + T^{6+2\rho-8\s\over\rho+1})\ll_\e T^{1+\e}
\cr}\leqno(6.4)
$$
for $\s \ge (5+\rho)/8$ if we choose $X = T^{4/(1+\rho)}$.
However, for $\rho\le 181/104 = 1.74038\ldots\,$ we have
$$
{1+\rho\over4} + {131\over832} \le {5+\rho\over8}.
$$
With $\rho=3/2$ we have $(5+\rho)/8 = 13/16$,
and   the first assertion of (2.14) follows. Note that from (see (2.12))
$$
\int_X^{2X}|\zx|^8x^{1-2\s}\d x \ll 1\qquad(\s > \hf(1+\rho)),
$$
one deduces easily that (see e.g.,  [9, Lemma 4])
$$
\int_T^{2T}|{\Cal Z}_2(\s + it)|^2\d t \ll 1
\qquad(\s > \hf(1+\rho)).\leqno(6.5)
$$
Thus the bound in (2.13) follows from (6.4), (6.5) and the convexity of mean
values for regular functions (cf. [6, Lemma 8.3]). Setting $\s=1, \rho = 3/2$ in
(2.13) we obtain the second bound in (2.14).
We remark that, by [9, eq. (3.21)], we have
$$
\int_T^{2T}|\zt|^8\d t \ll_\e T^{2\s-1}\int_0^{T^{1+\e}}|\Z_2(\s+it)|^2\d t
\quad(\hf < \s \le 1).\leqno(6.6)
$$
The bound (6.6) links the eighth moment of $|\zt|$ to the mean square of
$\Z_2(s)$. Thus the second bound in (2.14) gives the value $\rho = 11/7$,
which is close to the  best known bound $\rho = 3/2$.

%\vfill
%\eject
%\topskip2cm

\bigskip\bigskip
\Refs
\bigskip\bigskip

\item{[1]} F.V. Atkinson, The mean value of the Riemann zeta-function,
Acta Math. {\bf81}(1949), 353-376.

\item{[2]} J.B. Conrey, D.W. Farmer, J.P. Keating, M.O. Rubinstein
and N.C. Snaith, Integral moments of $L$-functions, Proc. Lond. Math.
Soc., III. Ser. {\bf91}(2005), 33-104.

\item{[4]} D.R. Heath-Brown, The twelfth power moment of the Riemann
zeta-function, Quart. J. Math. (Oxford) {\bf29}(1978), 443-462.

%% \item{[5]} M.N. Huxley, Exponential sums and the Riemann zeta-function V,
%% Proc. London Math. Soc. (3) {\bf90}(2005), 1-41.

%% \item{[6]} M.N. Huxley and A. Ivi\'c,
%% Subconvexity for the Riemann zeta-function and the divisor
%% problem, Bulletin  de l'Acad\'emie Serbe des Sciences et des
%% Arts - 2007,  Classe des Sciences math\'ematiques et naturelles,
%% Sciences math\'ematiques No. {\bf32}, pp. 13-32.

 \item{[5]}  A.E. Ingham, Mean-value theorems in the theory of the Riemann
 zeta-function,  Proc. London Math. Soc.  (2){\bf27}(1926), 273-300.

\item{[6]} A. Ivi\'c, The Riemann zeta-function, John Wiley \&
Sons, New York, 1985 (2nd ed., Dover, Mineola, N.Y., 2003).

\item{[7]} A. Ivi\'c, The mean values of the Riemann zeta-function,
LNs {\bf 82}, Tata Inst. of Fundamental Research, Bombay (distr. by
Springer Verlag, Berlin etc.), 1991.

\item{[8]} A. Ivi\'c,  On the fourth moment of the Riemann
zeta-function, Publs. Inst. Math. (Belgrade) {\bf 57(71)}(1995), 101-110.

\item{[9]} A. Ivi\'c, On some conjectures and
results for the Riemann zeta-function
and Hecke series, Acta Arith. {\bf109}(2001), 115-145.

\item{[10]}A. Ivi\'c, On the Riemann zeta function and the divisor problem,
Central European Journal of Mathematics 2({\bf4}) (2004),   1-15;
II ibid. 3({\bf2}) (2005), 203-214,
III, Annales Univ. Sci. Budapest., Sect. Comp. {\bf29}(2008), 3-23., and IV,
Uniform Distribution Theory {\bf1}(2006), 125-135.

\item{[11]} A. Ivi\'c, On moments of $|\zt|$ in short intervals,
Ramanujan Math. Soc. LNS{\bf2},
The Riemann zeta function and related themes: Papers in honour of
Professor Ramachandra, 2006, 81-97.

\item{[12]} A. Ivi\'c, On the estimation of some
Mellin transforms connected with the fourth
moment of $|\zeta({1\over2}+it)|$, Elementare und Analytische Zahlentheorie
(Tagungsband), Proceedings ELAZ-Conference May 24-28, 2004 (W. Schwarz
und J. Steuding eds.), Franz Steiner Verlag 2006, pp. 77-88.

\item{[13]} A. Ivi\'c, Some remarks on the moments of $|\zt|$ in short intervals,
Acta Math. Hung. {\bf119}(2008), 15-24.

\item{[14]} A. Ivi\'c, On the moments of the function $E^*(t)$,
Turkish Journal of Analysis and Number Theory, 2014, to appear.

\item{[15]} A. Ivi\'c, On the divisor function and the Riemann
zeta-function in short intervals, The Ramanujan Journal {\bf19}(2009), 207-224.

\item{[16]} A. Ivi\'c and Y. Motohashi, The mean square of the
error term for the fourth moment of the zeta-function,  Proc. London Math.
Soc. (3){\bf66}(1994), 309-329.

\item {[17]}  A. Ivi\'c and Y. Motohashi,  The fourth moment of the
Riemann zeta-function,  J. Number Theory  {\bf51}(1995), 16-45.

\item{[18]} A. Ivi\'c and P. Sargos, On the higher moments of the
error term in the divisor problem, Illinois J. Math. {\bf51}(2007),
353-377.

\item{[19]} A. Ivi\'c, M. Jutila  and Y. Motohashi,
The Mellin transform of powers of the zeta-function, Acta
Arith. {\bf95}(2000), 305-342.

\item{[20]} M. Jutila, On the divisor problem for short intervals,
Ann. Univer. Turkuensis Ser. {\bf A}I {\bf186}(1984), 23-30.

\item{[21]} M. Jutila, Riemann's zeta-function and the divisor problem,
Arkiv Mat. {\bf21}(1983), 75-96 and II, ibid. {\bf31}(1993), 61-70.

\item {[22]}  Y. Motohashi,  Spectral theory of the Riemann
zeta-function,  Cambridge University Press, Cambridge, 1997.

\item{[23]} K. Ramachandra, On the mean-value and omega-theorems
for the Riemann zeta-function, LNs {\bf85}, Tata Institute of Fundamental
Research, Bombay (distr. by Springer Verlag, Berlin etc.) 1995.

\item{[24]} O. Robert and P. Sargos, Three-dimensional exponential
sums with monomials, J. reine angew. Math. {\bf591}(2006), 1-20.

\item{[25]} E.C. Titchmarsh, The theory of the Riemann
zeta-function (2nd edition),  University Press, Oxford, 1986.

\item{[26]} N. Watt, A note on the mean square of $|\zt|$,
 J. London Math. Soc. {\bf82}(2)(2010), 279-294.

%\vskip1cm
\endRefs

\enddocument

\bye